\documentclass[a4paper, 12pt]{article}

\usepackage{amsmath}
\usepackage{amstext}
\usepackage{amsopn}
\usepackage{amssymb}
\usepackage{amsfonts}
\usepackage{amsthm}

\usepackage{graphicx}
\usepackage{fullpage}
\usepackage{url}
\usepackage{boxedminipage}
\usepackage[normalem]{ulem}
\usepackage{setspace}
\usepackage{enumerate}

\newcommand{\F}{\mathbb{F}}

\newcommand{\N}{\mathbf{N}}

\newcommand{\predfloat}{\operatorname{pred}}
\newcommand{\succfloat}{\operatorname{succ}}
\newcommand{\cond}{\operatorname{cond}}
\newcommand{\fl}{\operatorname{fl}}

\newcommand{\uu}{\ensuremath{{\mathbf{u}}}}

\newcommand{\eps}{\varepsilon}

\newcommand{\wh}[1]{{\,\widehat{#1}}}
\newcommand{\wt}[1]{\ensuremath{\,\widetilde{#1}}}
	  \newcommand{\wb}[1]{\ensuremath{\,\overline{#1}}}

\newcommand{\ie}{{\it i.e.},~}
\newcommand{\eg}{{\it e.g.},~}
\newcommand{\apriori}{{\em a priori} \,}

\theoremstyle{plain}
\newtheorem{theorem}{Theorem}

\newtheorem{lemma}[theorem]{Lemma}

\theoremstyle{remark}
\newtheorem{remark}{Remark}
\newtheorem{algorithm1}{Algorithm}
\newenvironment{algorithm}{\begin{algorithm1}\rm\singlespace}{\end{algorithm1}}

\newcommand{\Theorem}[1]{Theorem~\ref{#1}}
\newcommand{\Prop}[1]{Proposition~\ref{#1}}
\newcommand{\Lemma}[1]{Lemma~\ref{#1}}

\newcommand{\Algo}[1]{Algorithm~\ref{#1}}

\newcommand{\algo}[1]{{\sf #1}}
\newcommand{\apply}[2]{\mbox{\algo{#1}}\left(#2\right)}

\newcommand{\Sec}[1]{Section~\ref{#1}}

\newcommand{\Rel}[1]{Relation~(\ref{#1})}

\newcommand{\ineq}[1]{inequality~(\ref{#1})}
\newcommand{\Tab}[1]{Table~\ref{#1}}

\newcommand{\Fig}[1]{Figure~\ref{#1}}

\title{Faithful Polynomial Evaluation \\ with Compensated Horner Algorithm}
\date{\today}
\author{Philippe Langlois, Nicolas Louvet\\
  {\small Universit\'{e} de Perpignan Via Domitia%
    \thanks{DALI Research Team. Laboratory LP2A. 52, avenue Paul
      Alduy. F-66860 Perpignan, France.}}\\ 
{\small \textsf{\{langlois, nlouvet\}@univ-perp.fr}}
}

\begin{document}
%\doublespacing

\maketitle
\vspace{-1cm}
\begin{abstract}
This paper presents two sufficient conditions to ensure a faithful
evaluation of polynomial in IEEE-754 floating point
arithmetic. Faithfulness means that the computed value is one of the 
two floating point neighbours of the exact result; it can be satisfied
using a more accurate algorithm than the classic Horner scheme.
One condition here provided is an \apriori bound of the polynomial
condition number derived from the error analysis of the compensated
Horner algorithm. The second condition is both dynamic and validated
to check at the running time the faithfulness of a given
evaluation. Numerical experiments illustrate the behavior of these two
conditions and that associated running time over-cost is really
interesting.    
\end{abstract}


%%%%%%%%%%%%%%%%%%%%%%%%%%%%%%%%%%%%%%%%%%%%%%%%%%%%%%%%%%%%%%%%%%%%%%%%%%%%%%%%%%
\section{Introduction}
\label{sec:Introduction}
%%%%%%%%%%%%%%%%%%%%%%%%%%%%%%%%%%%%%%%%%%%%%%%%%%%%%%%%%%%%%%%%%%%%%%%%%%%%%%%%%%
\subsection{Motivation}
Horner's rule is the classic algorithm when evaluating a polynomial
$p(x)$. When performed in floating point arithmetic this algorithm
may suffer from (catastrophic) cancellations and so yields a computed value
with less exact digits than expected. The relative accuracy of the
computed value $\wh{p}(x)$ verifies the well known following inequality, 
\begin{equation}\label{rel:HornerBound}
\frac{|p(x) - \wh{p}(x)|}{|p(x)|} \leq \alpha(n) \cond(p, x) \  \uu.
\end{equation}  
In the right-hand side of this accuracy bound, $\uu$ is the computing precision
and $\alpha(n) \approx 2n$ for a polynomial of degree $n$. The
condition number $\cond(p,x)$ that only depends on $x$ and on $p$
coefficients will be explicited further. 
The product $\alpha(n) \cond(p, x)$ may be arbitrarily larger than
$1/\uu$ when cancellations appear, \ie when evaluating the polynomial
$p$ at the $x$ entry is ill-conditioned. \\

When the computing precision $\uu$ is not sufficient to guarantee a
desired accuracy, several solutions simulating a computation with
more bits exist. Priest-like ``double-double'' algorithms are
well-known and well-used solutions to simulate twice the IEEE-754 double 
precision \cite{Prie:91,XBLAS:02}. The compensated Horner algorithm
is a fast alternative to ``double-double'' introduced in
\cite{grll05} --- fast means that the compensated algorithm 
should run at least twice as fast as the ``double-double'' counterpart
with the same output accuracy. 
In both cases this accuracy is improved and now verifies 
\begin{equation}\label{rel:CompHornerBound}
\frac{|p(x)-\wh{p}(x)|}{|p(x)|} \leq \uu + \beta(n) \cond(p, x) \ \uu^2,
\end{equation}
with $\beta(n) \approx 4n^2$. This relation means that the computed
value is as accurate as the result of the Horner algorithm performed
in twice the working precision and then rounded to this working
precision.  \\ 
     
This bound also tells us that such algorithms may yield a full
precision accuracy for not too ill-conditioned polynomials, \eg when
$\beta(n) \cond(p, x) \uu < 1$. \\

This remark motivates this paper where we consider {\em faithful polynomial
  evaluation}. By faithful (rounding) we mean that the computed result
$\wh{p}(x)$ is one of the two floating point neighbours of the exact
result $p(x)$. Faithful rounding is known to be an interesting
property since for example it guarantees the correct sign
determination of arithmetic expressions, \eg for geometric
predicates.  

We first provide an \apriori sufficient criterion on the
condition number of the polynomial evaluation to ensure that the
compensated Horner algorithm provides a faithful rounding of the exact
evaluation (\Theorem{prop:FaithfulRoundingCompHorner} in
\Sec{sec:AprioriBound}).  
We also propose a validated and dynamic bound to prove
at the running time that the computed evaluation is actually
faithful (\Theorem{prop:DynErrorBounds} in \Sec{sec:DynBound}). 
We present numerical experiments to show that the dynamic
bound is sharper than the \apriori condition and we measure that the
corresponding over-cost is reasonable (\Sec{sec:ExperimentalResults}).

%%%%%%%%%%%%%%%%%%%%%%%%%%%%%%%%%%
\subsection{Notations}
%%%%%%%%%%%%%%%%%%%%%%%%%%%%%%%%%%
Throughout the paper, we assume a floating point arithmetic adhering to the
IEEE-754 floating point standard \cite{IEEE:85}. We constraint all the
computations to be performed in one working precision, with the ``round to the
nearest'' rounding mode. We also assume that no overflow nor underflow occurs
during the computations. Next notations are standard (see
\cite[chap.~2]{ASNA:02} for example). $\F$ is the set of all normalized floating
point numbers and $\uu$ denotes the unit roundoff, that is half the spacing
between $1$ and the next representable floating point value. For IEEE-754
double precision with rounding to the nearest, we have $\uu = 2^{-53} \approx 1.11
\cdot 10^{-16}$.
We define the floating point predecessor and successor of a real number $r$
as follows,
\begin{equation*}
\predfloat(r) = \max\{f\in\F / f<r\} \quad \mbox{and} \quad
\succfloat(r) = \min\{f\in\F / r<f\}.
\end{equation*}
A floating point number $f$ is defined to be a faithful rounding of a real
number $r$ if
\begin{equation*}
\predfloat(f) < r < \succfloat(f).
\end{equation*}

The symbols $\oplus$, $\ominus$, $\otimes$
and $\oslash$ represent respectively the floating point addition,
subtraction, multiplication and division. For more complex arithmetic
expressions, $\fl(\cdot)$ denotes the result of a floating point
computation where every operation inside the parenthesis is performed in
the working precision. So we have for example, $a \oplus b = \fl(a+b)$.

When no underflow nor overflow occurs, the following standard model describes
the accuracy of every considered floating point computation. For two floating
point numbers $a$ and $b$ and for $\circ$ in $\{+, -, \times, /\}$, the
floating point evaluation $\fl(a \circ b)$ of $a \circ b$ is such that
\begin{equation} \label{rel:StdModel}
\fl(a \circ b) = (a \circ b)(1+\eps_1) = (a \circ b)/(1+\eps_2),
\mbox{with} \quad |\eps_1|,|\eps_2| \leq \uu.
\end{equation}

To keep track of the $(1+\eps)$ factors in next error analysis, we use the
classic $(1+\theta_k)$ and $\gamma_k$ notations \cite[chap.~3]{ASNA:02}.
For any positive integer $k$, $\theta_k$ denotes a quantity bounded according to
\begin{equation*}
|\theta_k| \leq \gamma_k = \frac{k\uu}{1-k\uu}.
\end{equation*}
When using these notations, we always implicitly assume $k\uu<1$. In
further error analysis, we essentially use the following relations,
\begin{equation*}
(1+\theta_k)(1+\theta_j) \leq (1+\theta_{k+j}), \quad
k\uu \leq \gamma_k, \quad
\gamma_k \leq \gamma_{k+1}.
\end{equation*}
Next bounds are computable floating point values that will be useful
to derive dynamic validation in \Sec{sec:DynBound}.
We denotes $\fl(\gamma_k) = (k \uu) \oslash (1 \ominus k \uu)$ by $\wh{\gamma}_k$.
We know that $\fl(k\uu) = k\uu\in \F$, and $k\uu<1$ implies $\fl(1-k\uu) = 1-k\uu \in \F$.
So $\wh{\gamma}_k$ only suffers from a rounding error in the division and
\begin{equation} \label{rel:GammaK}
\gamma_k \leq (1+\uu) \wh{\gamma}_k.
\end{equation}
The next bound comes from the direct application of
Relation~(\ref{rel:StdModel}). For $x \in \F$ and $n \in \N$, 
\begin{equation} \label{rel:CompBound}
(1+\uu)^n|x| \leq
\fl \left(
\frac{|x|}{1-(n+1)\uu}
\right).
\end{equation}

%%%%%%%%%%%%%%%%%%%%%%%%%%%%%%%%%%%%%%%%%%%%%%%%%%%%%%%%%%%%%%%%%%%%%%%%%%%%%%%%%%
\section{From Horner to compensated Horner algorithm}
\label{sec:CompHorner}
%%%%%%%%%%%%%%%%%%%%%%%%%%%%%%%%%%%%%%%%%%%%%%%%%%%%%%%%%%%%%%%%%%%%%%%%%%%%%%%%%%
The compensated Horner algorithm improves the classic Horner iteration
computing a correcting term to compensate the rounding errors the
classic Horner iteration generates in floating point arithmetic. Main
results about compensated Horner algorithm are summarized in this
section; see \cite{grll05} for a complete description.  

%%%%%%%%%%%%%%%%%%%%%%%%%%%%%%%%%%
\subsection{Polynomial evaluation and Horner algorithm}
%%%%%%%%%%%%%%%%%%%%%%%%%%%%%%%%%%
The classic condition number of the evaluation of $p(x) = \sum_{i=0}^n a_i x^i$
at a given data $x$ is
\begin{equation} \label{rel:CondPoly}
\cond(p,x) = 
\frac{\sum_{i=0}^n |a_i| |x|^i}{|\sum_{i=0}^n a_i x^i|}=
\frac{\wt{p}(x)}{|p(x)|}.
\end{equation}
For any floating point value $x$ we denote by $\apply{Horner}{p, x}$
 the result of the floating point evaluation of the polynomial $p$ at
$x$ using next classic Horner algorithm. 

\begin{algorithm} \label{algo:Horner}
Horner algorithm
\begin{tabbing}
\quad \=\quad  \kill
function  $r_0 = \apply{Horner}{p,x}$ \\
$r_n = a_n$ \\
for $i=n-1:-1:0$ \\
\> $r_i = r_{i+1} \otimes x \oplus a_i$ \\
end
\end{tabbing}
\end{algorithm}

The accuracy of the result of \Algo{algo:Horner} verifies introductory
\ineq{rel:HornerBound} with $\alpha_n \uu= \gamma_{2n}$ and previous condition
number~(\ref{rel:CondPoly}). 
%% is linked to the condition number
%% of the polynomial evaluation thanks to the following forward error bound
%% (see~\cite[p.95]{ASNA:02}),
%% \begin{equation} \label{rel:REBHorner}
%% \frac{|p(x) - \apply{Horner}{p,x}|}{|p(x)|} \leq \gamma_{2n} \cond(p, x).
%% \end{equation}
Clearly, the condition number $\cond(p, x)$ can be arbitrarily large. In particular,
when $\cond(p, x) > 1/ \gamma_{2n}$, we cannot guarantee that the computed result
$\apply{Horner}{p,x}$ contains any correct digit.

We further prove that the error generated by the Horner algorithm is
exactly the sum of two polynomials with floating point coefficients. The next
lemma gives bounds of the generated error when evaluating this sum of
polynomials applying the Horner algorithm.

\begin{lemma} \label{lemma:HornerSum}
Let $p$ and $q$ be two polynomials with floating point coefficients,
such that $p(x) = \sum_{i=0}^n a_i x^i$ and $q(x) = \sum_{i=0}^n b_i x^i$.
We consider the floating point evaluation of $(p+q)(x)$ computed with
$\apply{Horner}{p \oplus q, x}$. Then, in case no
underflow occurs, the computed result satisfies the following forward
error bound,
%
%\begin{equation} \label{rel:ErrorBoundHornerSum}
%|(p+q)(x) - \apply{Horner}{p \oplus q, x}|
%\leq \gamma_{2n+1} (\wt{p}+\wt{q})(x).
%\end{equation}
%
\begin{equation} \label{rel:ErrorBoundHornerSum}
|(p+q)(x) - \apply{Horner}{p \oplus q, x}|
\leq \gamma_{2n+1} (\wt{p+q})(x).
\end{equation}
Moreover, if we assume that $x$ and the coefficients of $p$ and $q$ are
non-negative floating point numbers  then
\begin{equation} \label{rel:ComputableBoundHornerSum}
(p+q)(x) \leq
(1+\uu)^{2n+1} \apply{Horner}{p \oplus q, x}.
\end{equation}
\end{lemma}

\begin{proof}
The proof of the error bound~(\ref{rel:ErrorBoundHornerSum}) is easily adapted
from the one of the Horner algorithm (see \cite[p.95]{ASNA:02} for example).
To prove~(\ref{rel:ComputableBoundHornerSum}) we consider \Algo{algo:Horner},
where
\begin{equation*}
r_{n} = a_n \oplus b_n
\quad \mbox{and} \quad
r_{i} = r_{i+1} \otimes x \oplus (a_i \oplus b_i)
\quad \mbox{for} \quad
i=n-1,\ldots,0.
\end{equation*}
Next, using the standard model~(\ref{rel:StdModel}) it is easily proved by
induction that, for $i=0,\ldots,n$,
\begin{equation} \label{rel:ComputableBoundHornerSum:rel1}
\sum_{j=0}^{i}(a_{n-i+j}+b_{n-i+j})x^j \leq (1+\uu)^{2i+1} r_{n-i},
\end{equation}
which in turn proves~(\ref{rel:ComputableBoundHornerSum}) for $i=n$.
%Let us prove by induction that, for $i=0,\ldots,n$,
%\begin{equation} \label{rel:ComputableBoundHornerSum:rel1}
%\sum_{j=0}^{i}(a_{n-i+j}+b_{n-i+j})x^j \leq (1+\uu)^{2i+1} r_{n-i}.
%\end{equation}
%For $i=0$, $(a_n+b_n) \leq (1+\uu)(a_n \oplus b_n) = (1 + \uu) r_n$, so
%\Rel{rel:ComputableBoundHornerSum:rel1} is satisfied. Now we assume that
%\Rel{rel:ComputableBoundHornerSum:rel1} is true for some integer $i$ with
%$0 \leq i < n$. Then we have
%\begin{eqnarray*}
%\sum_{j=0}^{i+1}(a_{n-(i+1)+j}+b_{n-(i+1)+j})x^j
%& = &
%\sum_{j=0}^{i}(a_{n-i+j}+b_{n-i+j})x^j x + (a_{n-i}+b_{n-i})\\
%&\leq&
%(1+\uu)^{2i+1} r_{n-i} x + (a_{n-i}+b_{n-i})\\
%&\leq&
%(1+\uu)^{2i+1} (1+\uu)^2
%\left( r_{n-i} \otimes x \oplus (a_{n-i} \oplus b_{n-i}) \right)\\
%&\leq&
%(1+\uu)^{2(i+1)+1}
%r_{n-(i+1)}.
%\end{eqnarray*}
%Therefore \Rel{rel:ComputableBoundHornerSum:rel1} is proved by induction, which
%in turn proves~(\ref{rel:ComputableBoundHornerSum}).
\end{proof}

%%%%%%%%%%%%%%%%%%%%%%%%%%%%%%%%%%%%%%%%%%%%%%%%%%%%%%%%%%%%%%%%%%%%%%%%%%%%%%%%%%
%\section{Error free transformations} \label{sec:EFT}
%%%%%%%%%%%%%%%%%%%%%%%%%%%%%%%%%%%%%%%%%%%%%%%%%%%%%%%%%%%%%%%%%%%%%%%%%%%%%%%%%%

%%%%%%%%%%%%%%%%%%%%%%%%%%%%%%%%%%
\subsection{EFT for the elementary operations}
%%%%%%%%%%%%%%%%%%%%%%%%%%%%%%%%%%
Now we review well known results concerning error free
transformation (EFT) of the elementary floating point operations $+$, $-$ and
$\times$.\\

% We also introduce an EFT for the polynomial evaluation with the
% Horner algorithm.\\ 

Let $\circ$ be an operator in $\{ +, -, \times\}$, $a$ and $b$ be two floating
point numbers, and $\wh{x} = \fl(a \circ b)$. Then their exist a floating point
value $y$ such that
\begin{equation} \label{rel:DefEFT}
a \circ b = \wh{x} + y.
\end{equation}
The difference $y$ between the exact result and the computed result is 
the rounding error generated by the computation of $\wh{x}$.  Let us emphasize
that relation~(\ref{rel:DefEFT}) between four floating point values relies on
real operators and exact equality, \ie not on approximate floating point
counterparts. Ogita {\it et al.}~\cite{OgRO:05} name such a transformation an
error free transformation (EFT). The practical interest of the EFT comes from
next Algorithms \ref{algo:TwoSum} and \ref{algo:TwoProd} that compute
the exact error term $y$ for $\circ = +$ 
and $\circ = \times$.\\
%Let $\circ$ be in $\{ +, -, \times\}$, $a$ and $b$ be two floating
%point numbers, and $\wh{x} = \fl(a \circ b)$. The {\it elementary rounding error}
%in the computation of $\wh{x}$ is
%\begin{equation}
%y = (a \circ b) - \fl(a \circ b),
%\end{equation}
%that is the difference between the exact result and the computed result of
%the operation. In particular, for $\circ$ in $\{ +, -, \times\}$, the elementary
%rounding error $y$ both belongs to $\F$, and is computable using only the operations
%defined within $\F$. Thus, for $\circ$ in $\{ +, -, \times\}$,
%any pair of inputs $(a, b)$ in $\F^2$ can be transformed into an output pair
%$(\wh{x}, y)$ in $\F^2$ such that
%\begin{displaymath}
%a \circ b = \wh{x} + y \quad \mbox{and} \quad \wh{x} = \fl(a \circ b).
%\end{displaymath}
%Let us emphasize that this relation between these four floating point values
%relies on real operators and exact equality (\ie not on approximate floating
%point counterparts). Ogita {\it et al.}~\cite{OgRO:05} call such a transformation
%an {\it error free transformation} (EFT).
%

For the EFT of the addition we use \Algo{algo:TwoSum}, the well known
\algo{TwoSum} algorithm by Knuth~\cite{Knut:98} that requires 6 flop (floating
point operations).
For the EFT of the product, we first need to split the input arguments into two
parts. It is done using \Algo{algo:Split} of Dekker~\cite{Dekk:71}
where $r=27$ for IEEE-754 double precision. Next, \Algo{algo:TwoProd} by Veltkamp (see~\cite{Dekk:71}) can be
used for the EFT of the product. This algorithm is commonly called
\algo{TwoProd} and requires 17 flop. 

\begin{center}
\begin{minipage}[c]{0.8\linewidth}
\vspace{0.5cm}\hbox{\raisebox{0.2em}{\vrule depth 0pt height 0.4pt width \linewidth}}
\begin{algorithm}
EFT of the sum of two floating point numbers. \label{algo:TwoSum}\vspace{-0.2cm}
\begin{tabbing}
\quad \=\quad \=\quad \kill
function $[x, y] = \apply{TwoSum}{a, b}$\\
\> $x = a \oplus b$\\
\> $z = x \ominus a$\\
\> $y = (a \ominus (x \ominus z)) \oplus (b \ominus z)$
\end{tabbing}
\end{algorithm}
\vspace{-1em}\hbox{\raisebox{0.3em}{\vrule depth 0pt height 0.4pt width \linewidth}}\vspace{-1em}
\begin{algorithm}
Splitting of a floating point number into two parts. \label{algo:Split}\vspace{-0.2cm}
\begin{tabbing}
\quad \=\quad \=\quad \kill
function $[x, y] = \apply{Split}{a}$\\
\> $z = a \otimes (2^r + 1)$\\
\> $x = z \ominus (z \ominus a)$\\
\> $y = a \ominus x$
\end{tabbing}
\end{algorithm}
\vspace{-1em}\hbox{\raisebox{0.3em}{\vrule depth 0pt height 0.4pt width \linewidth}}\vspace{-1em}
\begin{algorithm}
EFT of the product of two floating point numbers. \label{algo:TwoProd}\vspace{-0.2cm}
\begin{tabbing}
\quad \=\quad \=\quad \kill
function $[x, y] = \apply{TwoProd}{a, b}$\\
\> $x = a \otimes b$\\
\> $[a_h, a_l] = \apply{Split}{a}$\\
\> $[b_h, b_l] = \apply{Split}{b}$\\
\> $y = a_l \otimes b_l \ominus (((x \ominus a_h \otimes b_h) \ominus a_l \otimes b_h) \ominus a_h \otimes b_l )$
\end{tabbing}
\end{algorithm}
\vspace{-1em}\hbox{\raisebox{0.3em}{\vrule depth 0pt height 0.4pt
    width \linewidth}}%\vspace{-1em} 
%
%% \begin{algorithm}
%% EFT of the product of two floating point numbers with a FMA. \label{algo:TwoProdFMA}\vspace{-0.2cm}
%% \begin{tabbing}
%% \quad \=\quad \=\quad \kill
%% function $[x, y] = \apply{TwoProdFMA}{a, b}$\\
%% \> $x = a \otimes b$\\
%% \> $y = \apply{FMA}{a, b, -x}$
%% \end{tabbing}
%% \end{algorithm}
%% %
%% \vspace{-1em}\hbox{\raisebox{0.3em}{\vrule depth 0pt height 0.4pt width \linewidth}}\vspace{0.5cm}
\end{minipage}
\end{center}

%% If $q$ is the number of bits of the mantissa, let $r = \lceil q/2 \rceil$.
%% \Algo{algo:Split} splits a floating point number $a$ into two parts $x$ and $y$,
%% both having at most $r-1$ nonzero bits, such that $a = x + y$. 
%% For example, with the IEEE-754 double precision, $q = 53$, $r = 27$, therefore the output numbers
%% have at most $26$ bits. The trick is that one bit sign is used for the 
%% splitting. 

The next theorem exhibits the previously announced properties of \algo{TwoSum}
and \algo{TwoProd}.
\begin{theorem}[\cite{OgRO:05}] \label{theorem:EFT}
Let $a, b$ in $\F$ and $x, y \in \F$ such that $[x, y] = \algo{TwoSum}(a,b)$
(Algorithm \ref{algo:TwoSum}). Then, ever in the presence of underflow,
\begin{equation*}
a + b = x + y,
\quad x = a \oplus b,
\quad |y| \leq \uu |x|,
\quad |y| \leq \uu |a+b|.
\end{equation*}
Let $a, b \in \F$ and $x, y \in \F$ such that $[x, y] = \algo{TwoProd}(a,b)$
(Algorithm \ref{algo:TwoProd}). Then, if no underflow occurs,
\begin{equation*}
a \times b = x + y,
\quad x = a \otimes b,
\quad |y| \leq \uu |x|,
\quad |y| \leq \uu |a \times b|.
\end{equation*}
\end{theorem}

%% \algo{TwoProd} can be rewritten straightforwardly for  processors that provide a
%% Fused-Multiply-and-Add operator (\algo{FMA}), such as Intel Itanium or IBM
%% PowerPC~\cite{Niev:03,OgRO:05}. For $a$, $b$ and $c$ in $\F$,
%% $\apply{FMA}{a,b,c}$ is the exact result $a \times b + c$ rounded to the nearest
%% floating point value. Thus $y=a \times b - a \otimes b = \apply{FMA}{a, b, -(a
%% \otimes b)}$ and \algo{TwoProd} can be replaced by \Algo{algo:TwoProdFMA},
%% requiring now only 2 flop.

We notice that algorithms \algo{TwoSum} and \algo{TwoProd} 
%and \algo{TwoProdFMA}
only require well optimizable floating point operations. They do not use branches,
nor access to the mantissa that can be time-consuming.
We just mention that significant improvements of these algorithms are
defined when a Fused-Multiply-and-Add operator is available \cite{grll05}. 

%%%%%%%%%%%%%%%%%%%%%%%%%%%%%%%%%%%%%%%%%%%%%%%%%%%%%%%%%%%%%%%%%%%%%%%%%%%%%%%
\subsection{An EFT for the Horner algorithm}
%%%%%%%%%%%%%%%%%%%%%%%%%%%%%%%%%%%%%%%%%%%%%%%%%%%%%%%%%%%%%%%%%%%%%%%%%%%%%%%
As previously mentioned, next EFT for the polynomial evaluation with
the Horner algorithm exhibits the exact rounding error generated by the
Horner algorithm together with an algorithm to compute it.

\begin{algorithm} \label{algo:EFTHorner}
\rm
EFT for the Horner algorithm
\begin{tabbing}
\quad \=\quad  \kill
function $[s_0, p_{\pi}, p_{\sigma}] = \algo{EFTHorner}(p,x)$ \\
$s_n = a_n$ \\
for $i=n-1:-1:0$ \\
\> $[p_i , \pi_i ]$ = \texttt{TwoProd}$(s_{i+1},x)$ \\
\> $[s_i , \sigma_i ]$ = \texttt{TwoSum}$(p_i,a_i)$ \\
\> Let $\pi_i$ be the coefficient of degree $i$ in $p_{\pi}$ \\
\> Let $\sigma_i$ be the coefficient of degree $i$ in $p_{\sigma}$ \\
end
\end{tabbing}
\end{algorithm}

\begin{theorem}[\cite{grll05}] \label{theorem:EFTHorner}
Let $p(x) = \sum_{i=0}^{n} a_i x^i$ be a polynomial of degree $n$ with
floating point coefficients, and let $x$ be a floating point value. Then
\Algo{algo:EFTHorner} computes both
\begin{enumerate}[i)]
\item the floating point evaluation $\apply{Horner}{p, x}$ and
\item two polynomials $p_{\pi}$ and $p_{\sigma}$ of degree $n-1$ with
floating point coefficients,
\end{enumerate}
such that
\begin{equation*}
[\apply{Horner}{p, x}, p_{\pi}, p_{\sigma}] = \apply{EFTHorner}{p, x}.
\end{equation*}
If no underflow occurs,
\begin{equation} \label{rel:EFTHorner}
p(x) = \apply{Horner}{p, x} + (p_{\pi}+p_{\sigma})(x).
\end{equation}
Moreover,
%
%\begin{equation} \label{rel:BoundErrorTerm}
%(\wt{p_{\pi}}+\wt{p_{\sigma}})(x) \leq \gamma_{2n} \wt{p}(x).
%\end{equation}
%
\begin{equation} \label{rel:BoundErrorTerm}
(\wt{p_{\pi}+p_{\sigma}})(x) \leq \gamma_{2n} \wt{p}(x).
\end{equation}
\end{theorem}

\Rel{rel:EFTHorner} means that algorithm \algo{EFTHorner} is an EFT
for polynomial evaluation with the Horner algorithm.

\begin{proof}[Proof of \Theorem{theorem:EFTHorner}]
Since \algo{TwoProd} and \algo{TwoSum} are EFT from \Theorem{theorem:EFT}
it follows that $s_{i+1} x = p_i + \pi_i$ and $p_i + a_i = s_i + \sigma_i$.
Thus we have $s_i = s_{i+1} x + a_i - \pi_i - \sigma_i$, for $i=0,\ldots,n-1$.
Since $s_n = a_n$, at the end of the loop we have
\begin{equation*}
s_0 = \sum_{i=0}^{n}a_i x^i-\sum_{i=0}^{n-1}\pi_i x^i-\sum_{i=0}^{n-1}\sigma_i x^i,
\end{equation*}
which proves~(\ref{rel:EFTHorner}).
%
%Now we prove relation~(\ref{rel:BoundErrorTerm}).
%From \Theorem{theorem:EFT}, since \algo{TwoSum} and \algo{TwoProd} are
%EFT, for $i=0,\cdots,n-1$, we have $|\pi_{i}| \leq \uu |p_{i}|$ and
%$|\sigma_{i}| \leq \uu |s_i|$. Therefore
%\begin{equation*}
%(\wt{p_{\pi}} + \wt{p_{\sigma}})(x)
%= \sum_{i=0}^{n-1}(|\pi_i|+|\sigma_i|)|x^i|
%\leq \uu \sum_{i=1}^{n}(|p_{n-i}|+|r_{n-i}|)|x^{n-i}|.
%\end{equation*}
%Applying the standard model of floating point arithmetic to the computations in the
%inner loop of \Algo{algo:EFTHorner}, it is easily proved by induction that
%\begin{equation*}
%|p_{n-i}| \leq (1+\gamma_{2i-1}) \sum_{j=1}^{i}|a_{n-i+j}| |x^j|
%\;\; \mbox{and} \;\;
%|r_{n-i}| \leq (1+\gamma_{2i}) \sum_{j=0}^{i}|a_{n-i+j}| |x^j|
%\;\; \mbox{for} \;\; i=1,\ldots,n.
%\end{equation*}
%Therefore we have $|p_{n-i}||x^i| \leq (1+\gamma_{2i-1}) \wt{p}(x)$ and
%$|r_{n-i}||x^i| \leq (1+\gamma_{2i}) \wt{p}(x)$ for $i=1,\ldots, n$. Hence
%\begin{equation*}
%\sum_{i=1}^{n}(|p_{n-i}|+|r_{n-i}|)|x^{n-i}|
%\leq \sum_{i=1}^{n}\left(2+\gamma_{2i-1}+\gamma_{2i}\right) \wt{p}(x)
%\leq 2 n \left(1+\gamma_{2n-2}\right) \wt{p}(x).
%\end{equation*}
%Thus
%$(\wt{p_{\pi}}+\wt{p_{\sigma}})(x) \leq 2n\uu\left(1+\gamma_{2n-2}\right)\wt{p}(x)$.
%Since $2n\uu(1+\gamma_{2n-2}) \leq \gamma_{2n}$, we finally
%obtain~(\ref{rel:BoundErrorTerm}).

Now we prove relation~(\ref{rel:BoundErrorTerm})
According to the error analysis of the Horner algorithm
(see \cite[p.95]{ASNA:02}), we can write
\begin{equation*}
\apply{Horner}{p, x} = (1+\theta_{2n})a_n x^n + \sum_{i=0}^{n-1} (1+\theta_{2i+1}) a_i x^i,
\end{equation*}
where every $\theta_{k}$ satisfies $|\theta_k| \leq \gamma_k$. 
Then using (\ref{rel:EFTHorner}) we have
\begin{equation*}
(p_{\pi} + p_{\sigma})(x)
  = p(x) - \apply{Horner}{p, x}
  = -\theta_{2n} a_n x^n - \sum_{i=0}^{n-1} \theta_{2i+1} a_i x^i.
\end{equation*}
Therefore it yields next expected inequalities between the absolute values, 
\begin{equation*}
(\wt{p_{\pi} + p_{\sigma}})(x) \leq \gamma_{2n} |a_n| |x|^n +
  \sum_{i=0}^{n-1} \gamma_{2i+1} |a_i| |x|^i 
  \leq \gamma_{2n}\wt{p}(x).
\end{equation*}
%and we finally obtain $(\wt{p_{\pi} + p_{\sigma}})(x) \leq \gamma_{2n}\wt{p}(x)$.
\end{proof}

%%%%%%%%%%%%%%%%%%%%%%%%%%%%%%%%%%%%%%%%%
\subsection{Compensated Horner algorithm}
%%%%%%%%%%%%%%%%%%%%%%%%%%%%%%%%%%%%%%%%%

From \Theorem{theorem:EFTHorner} the final forward error of the
floating point evaluation of $p$ at $x$ according to the Horner algorithm is
\begin{equation*}
c = p(x) - \apply{Horner}{p, x} = (p_{\pi} + p_{\sigma})(x),
\end{equation*}
where the two polynomials $p_{\pi}$ and $p_{\sigma}$ are exactly identified
by \algo{EFTHorner} (\Algo{algo:EFTHorner}) ---this latter also computes
$\apply{Horner}{p, x}$. 
Therefore, the key of the compensated algorithm is to
compute, in the working precision, first an approximate $\wh{c}$ of
the final error $c$ and then a corrected result 
\begin{equation*}
\wb{r} = \apply{Horner}{p, x} \oplus \wh{c}.
\end{equation*}
These two computations leads to next compensated Horner algorithm
\algo{CompHorner} (\Algo{algo:CompHorner}). 

\begin{algorithm} \label{algo:CompHorner}
Compensated Horner algorithm
\begin{tabbing}
\quad \=\quad  \kill
function $\wb{r} = \apply{CompHorner}{p, x}$ \\
$\left[\wh{r}, p_{\pi}, p_{\sigma}\right]
= \apply{EFTHorner}{p, x}$ \\
$\wh{c} = \apply{Horner}{p_{\pi} \oplus p_{\sigma}, x}$\\
$\wb{r} = \wh{r} \oplus \wh{c}$
\end{tabbing}
\end{algorithm}

We say that $\wh{c}$  is a correcting term for $\apply{Horner}{p, x}$.
The corrected result $\bar{r}$ is expected to be more accurate than the first
result $\apply{Horner}{p, x}$ as proved in next section.

%%%%%%%%%%%%%%%%%%%%%%%%%%%%%%%%%%%%%%%%%%%%%%%%%%%%%%%%%%%%%%%%%%%%%%%%%%%%%%%%%%
\section{An {\it a priori} condition for faithful rounding}
\label{sec:AprioriBound}

We start proving the accuracy behavior of the compensated Horner
algorithm we previously mentioned with introductory
inequality~(\ref{rel:CompHornerBound}) and that motivates the search
for a faithful polynomial evaluation. This bound (and its proof) is
the first step towards the proposed \apriori sufficient condition for
a faithful rounding with compensated Horner algorithm.

%%%%%%%%%%%%%%%%%%%%%%%%%%%%%%%%%%%%%%%%%%%%%%%%
\subsection{Accuracy of the compensated Horner algorithm}
%%%%%%%%%%%%%%%%%%%%%%%%%%%%%%%%%%%%%%%%%%%%%%%%
Next result proves that the result of a polynomial evaluation computed
with the compensated Horner algorithm (\Algo{algo:CompHorner}) is as accurate
as if computed by the classic Horner algorithm using twice the working precision
and then rounded to the working precision.

\begin{theorem}[\cite{grll05}] \label{theorem:CompHorner}
Consider a polynomial $p$ of degree $n$ with floating point coefficients,
and $x$ a floating point value. If no underflow occurs,
\begin{equation} \label{rel:EBCompHorner}
|\apply{CompHorner}{p, x} - p(x)| \leq \uu |p(x)| + \gamma_{2n}^{2} \wt{p}(x).
\end{equation}
\end{theorem}

\begin{proof}
The absolute forward error generated by \Algo{algo:CompHorner} is
\begin{equation*}
| \wb{r} - p(x) |
= \left| (\wh{r} \oplus \wh{c}) - p(x) \right|
= \left| (1 +\eps)(\wh{r} + \wh{c}) - p(x) \right| \quad \mbox{with} \quad |\eps| \leq \uu.
\end{equation*}
Let $c = (p_{\pi}+p_{\sigma})(x)$. From
\Theorem{theorem:EFTHorner} we have
$\wh{r} =\apply{Horner}{p, x} = p(x) - c$, thus
\begin{equation*}
| \wb{r} - p(x) |
=\left| (1 +\eps) \left( p(x) - c + \wh{c}\right) - p(x) \right|
\leq \uu |p(x)| + (1 + \uu) | \wh{c} - c |.
\end{equation*}
%
%Since $\wh{c} = \apply{Horner}{p_{\pi} \oplus p_{\sigma}, x}$ with
%$p_{\pi}$ and $p_{\sigma}$ two polynomials of degree $n-1$,
%\Lemma{lemma:HornerSum} yields
%$|\wh{c} - c| \leq \gamma_{2n-1} (\wt{p_{\pi}}+ \wt{p_{\sigma}})(x)$.
%Then we use \Theorem{theorem:EFTHorner} to bound $(\wt{p_{\pi}}+
%\wt{p_{\sigma}})(x)$ as
%\begin{equation} \label{rel:ProofCompHorner:rel1}
%|\wh{c} - c| \leq \gamma_{2n-1} \gamma_{2n} \wt{p}(x).
%\end{equation}
%Since $(1+\uu)\gamma_{2n-1} \leq \gamma_{2n}$, we finally write the expected
%error bound~(\ref{rel:EBCompHorner}).
%
Since $\wh{c} = \apply{Horner}{p_{\pi} \oplus p_{\sigma}, x}$ with
$p_{\pi}$ and $p_{\sigma}$ two polynomials of degree $n-1$,
\Lemma{lemma:HornerSum} yields
$|\wh{c} - c| \leq \gamma_{2n-1} (\wt{p_{\pi}+ p_{\sigma}})(x)$.
Then using (\ref{rel:BoundErrorTerm}) we have
$|\wh{c} - c| \leq \gamma_{2n-1} \gamma_{2n} \wt{p}(x)$.
Since $(1+\uu)\gamma_{2n-1} \leq \gamma_{2n}$, we finally write the expected
error bound~(\ref{rel:EBCompHorner}).
\end{proof}

\begin{remark}
For later use, we notice that
$|\wh{c} - c| \leq \gamma_{2n-1} \gamma_{2n} \wt{p}(x)$
implies
\begin{equation} \label{rel:EBCorrectingTerm}
|\wh{c} - c| \leq \gamma_{2n}^2 \wt{p}(x).
\end{equation}
\end{remark}

%%TODO: un peu une repetition de l'intro (avec des ecarts de
%%notation). Si besoin remplacer par :
% The previous theorem exhibits introductory
% \Rel{rel:CompHornerBound}. As mentionned before it describes two
% properties about the accuracy of compensated Horner. First this
% accuracy is similar to doubling the working precision for classic
% Horner algorithm. Therefore compensated Horner algorithm computes an
% evaluation accurate to the last few bits as long as the condition
% number is smaller than $\gamma_{2n}^{-1} \approx (2n\uu)^{-1}$. 

It is interesting to interpret the previous theorem in terms of the condition
number of the polynomial evaluation of $p$ at $x$. Combining the error bound
(\ref{rel:EBCompHorner}) with the condition number~(\ref{rel:CondPoly}) of
polynomial evaluation gives the precise writing of our introductory
\ineq{rel:CompHornerBound}, 
\begin{equation} \label{rel:REBCompHorner}
\frac{| \apply{CompHorner}{p, x} - p(x) |}{|p(x)|}
\leq \uu + \gamma_{2n}^{2} \cond(p, x).
\end{equation}
In other words, the bound for the relative error of the computed result is
essentially $\gamma_{2n}^2$ times the condition number of the polynomial
evaluation, plus the inevitable summand $\uu$ for rounding the result to the
working precision. In particular, if $\cond(p, x) < \uu/ \gamma_{2n}^2$, then the
relative accuracy of the result is bounded by a constant of the order $\uu$. This
means that the compensated Horner algorithm computes an evaluation accurate to the
last few bits as long as the condition number is smaller than $\uu/
\gamma_{2n}^2 \approx 1/4n^2\uu$. 
Besides that, relation~(\ref{rel:REBCompHorner}) tells us
that the computed result is as accurate as if computed by the classic Horner
algorithm with twice the working precision, and then rounded to the working
precision.

%%%%%%%%%%%%%%%%%%%%%%%%%%%%%%%%%%%%%%%%%%%%%%%%
\subsection{An \apriori condition for faithful rounding}
%%%%%%%%%%%%%%%%%%%%%%%%%%%%%%%%%%%%%%%%%%%%%%%%

Now we propose a sufficient condition on $\cond(p, x)$ to ensure that the
corrected result $\wb{r}$ computed with the compensated Horner algorithm is a
faithful rounding of the exact result $p(x)$. For this purpose, we use the
following lemma from \cite{RuOO06}.

\begin{lemma}[\cite{RuOO06}] \label{lemma:FaithfulRounding}
Let $r, \delta$ be two real numbers and $\wb{r} = \fl(r)$. We assume here that
$\wb{r}$ is a normalized floating point number. If
$|\delta| < \frac{\uu}{2} |\wh{r}|$ then $\wb{r}$ is a faithful rounding of
$r +\delta$.
\end{lemma}

From \Lemma{lemma:FaithfulRounding}, we derive a useful criterion to ensure that
the compensated result provided by \algo{CompHorner} is faithfully rounded to
the working precision.

\begin{lemma} \label{lemma:FaithfulRoundingCompHorner}
Let $p$ be a polynomial of degree $n$ with floating point coefficients,
and $x$ be a floating point value. We consider the approximate $\wb{r}$ of
$p(x)$ computed with $\apply{CompHorner}{p, x}$, and we assume that no
underflow occurs during the computation. Let $c$ denotes
$c = (p_{\pi} + p_{\sigma})(x)$. If $|\wh{c} - c| < \frac{\uu}{2}|\wb{r}|$,
then $\wb{r}$ is a faithful rounding of $p(x)$.
\end{lemma}

\begin{proof}
We assume that $|\wh{c} - c| < \frac{\uu}{2}|\wb{r}|$. From the notations
of~\Algo{algo:CompHorner}, we recall that $\fl(\wh{r} + \wh{c}) = \wb{r}$.
Then from \Lemma{lemma:FaithfulRounding} it follows that $\wb{r}$
is a faithful rounding of $\wh{r}+\wh{c}+c-\wh{c} = \wh{r}+c$.
Since $[\wh{r}, p_{\pi}, p_{\sigma}] = \apply{EFTHorner}{p, x}$,
\Theorem{theorem:EFTHorner} yields $p(x) = \wh{r}+c$.
Therefore $\wb{r}$ is a faithful rounding of $p(x)$.
\end{proof}

The criterion proposed in \Lemma{lemma:FaithfulRoundingCompHorner}
concerns the accuracy of the correcting term $\wh{c}$. Nevertheless
\Rel{rel:EBCorrectingTerm} pointed after the proof of 
\Theorem{theorem:CompHorner} says that the absolute 
error $|\wh{c} - c|$ is bounded by $\gamma_{2n}^2 \wt{p}(x)$.
This provides us a more useful criterion, since it relies on the condition number
$\cond(p, x)$, to ensure that \algo{CompHorner} computes a faithfully
rounded result.

\begin{theorem} \label{prop:FaithfulRoundingCompHorner}
Let $p$ be a polynomial of degree $n$ with floating point coefficients,
and $x$ a floating point value. If
\begin{equation} \label{rel:APrioriCriterion}
\cond(p, x) < \frac{1-\uu}{2+\uu} \uu {\gamma_{2n}}^{-2},
\end{equation}
then $\apply{CompHorner}{p, x}$ computes a faithful rounding of the
exact $p(x)$.
\end{theorem}

\begin{proof}
We assume that~(\ref{rel:APrioriCriterion}) is satisfied and we use the same
notations as in \Lemma{lemma:FaithfulRoundingCompHorner}.

First we notice that $\wb{r}$ and $p(x)$ are of the same sign. Indeed,
from~(\ref{rel:EBCompHorner}) it follows that
$\left|{\wb{r}}/{p(x)} - 1\right| \leq \uu + \gamma_{2n}^2 \cond(p, x)$,
and therefore ${\wb{r}}/{p(x)} \geq 1 - \uu -\gamma_{2n}^2 \cond(p, x)$.
But~(\ref{rel:APrioriCriterion}) implies that
$1 - \uu -\gamma_{2n}^2 \cond(p, x) > 1 - 3\uu/(2+\uu) > 0$,
hence ${\wb{r}}/{p(x)} > 0$.  Since $\wb{r}$ and $p(x)$ have the same sign,
it is easy to see that
\begin{equation} \label{proof:FaithfulRoundingCompHorner:rel1}
(1-\uu)|p(x)| - \gamma_{2n}^2 \wt{p}(x) \leq |\wb{r}|.
\end{equation}
Indeed, if $p(x)>0$ then~(\ref{rel:EBCompHorner}) implies
$p(x) - \uu |p(x)| - \gamma_{2n}^2 \wt{p}(x) \leq \wb{r} = |\wb{r}|$.
If $p(x)<0$, from~(\ref{rel:EBCompHorner})
it follows that $\wb{r} \leq p(x) + \uu |p(x)| + \gamma_{2n}^2 \wt{p}(x)$,
hence $- p(x) - \uu |p(x)| - \gamma_{2n}^2 \wt{p}(x) \leq - \wb{r} = |\wb{r}|$.

Next, a small computation proves that
\begin{equation*}
\cond(p, x) < \frac{1-\uu}{2+\uu} \uu {\gamma_{2n}}^{-2}
\quad \mbox{if and only if} \quad
\gamma_{2n}^2 \wt{p}(x)
< \frac{\uu}{2} \left[(1-\uu)|p(x)| - \gamma_{2n}^2 \wt{p}(x)\right].
\end{equation*}
Finally, from~(\ref{rel:EBCorrectingTerm})
and~(\ref{proof:FaithfulRoundingCompHorner:rel1}) it follows
\begin{equation*}
|\wh{c} - c|
\leq \gamma_{2n}^2 \wt{p}(x)
< \frac{\uu}{2} \left[ (1-\uu)|p(x)| - \gamma_{2n}^2 \wt{p}(x) \right]
\leq \frac{\uu}{2} |\wb{r}|.
\end{equation*}
From \Lemma{lemma:FaithfulRoundingCompHorner} we deduce that
$\wb{r}$ is a faithful rounding of $p(x)$.
\end{proof}

Numerical values of condition numbers for a faithful polynomial
evaluation in IEEE-754 double precision are presented in
\Tab{tab:BoundsCond} for degrees varying from 10 to 500.
\begin{table} 
\caption{{\it A priori} bounds on the condition number to ensure
  faithful rounding in IEEE-754 double precision for polynomials of degree
  10 to 500}
\label{tab:BoundsCond}
\begin{center}\vspace{-0.2cm}
\begin{tabular}{cccccccc}
n & 10 & 100 & 200 & 300 & 400 & 500\\
\hline
$\frac{1-\uu}{2-\uu} \uu {\gamma_{2n}}^{-2}$ &
$1.13 \cdot 10^{13}$ &
$1.13 \cdot 10^{11}$ &
$2.82 \cdot 10^{10}$ &
$1.13 \cdot 10^{10}$ &
$7.04 \cdot 10^{9}$  &
$4.51 \cdot 10^{9}$
\end{tabular}
\end{center}
\vspace{-0.5cm}
\end{table}

%%%%%%%%%%%%%%%%%%%%%%%%%%%%%%%%%%%%%%%%%%%%%%%%%%%%%%%%%%%%%%%%%%%%%%%%%%%%%%%%%%
\section{Dynamic and validated error bounds for faithful rounding and
  accuracy} 
\label{sec:DynBound}
%%%%%%%%%%%%%%%%%%%%%%%%%%%%%%%%%%%%%%%%%%%%%%%%%%%%%%%%%%%%%%%%%%%%%%%%%%%%%%%%%%

The results presented in \Sec{sec:AprioriBound} are perfectly suited for
theoretical purpose, for instance when we can \apriori bound the condition
number of the evaluation. However, neither the error bound in
\Theorem{theorem:CompHorner}, nor the criterion proposed in
\Theorem{prop:FaithfulRoundingCompHorner} can be easily checked using only
floating point arithmetic. Here we provide dynamic counterparts of
\Theorem{theorem:CompHorner} and \Prop{prop:FaithfulRoundingCompHorner},
that can be evaluated using floating point arithmetic in the ``round
to the nearest'' rounding mode.

\begin{lemma} \label{lemma:DynBoundErrorTerm}
Consider a polynomial $p$ of degree $n$ with floating point coefficients,
and $x$ a floating point value. We use the notations of
\Algo{algo:CompHorner}, and we denote $(p_{\pi}+p_{\sigma})(x)$ by $c$.
Then
\begin{equation} \label{rel:DynBoundErrorTerm}
|c-\wh{c}|
\leq \fl \left(
\frac{ \wh{\gamma}_{2n-1} \apply{Horner}{|p_{\pi}|\oplus|p_{\sigma}|,|x|} }{ 1-2(n+1)\uu }
\right)
:= \wh{\alpha}.
\end{equation}
\end{lemma}

\begin{proof}
Let us denote $\apply{Horner}{|p_{\pi}|\oplus|p_{\sigma}|,|x|}$ by $\wh{b}$.
Since $c = (p_{\pi}+p_{\sigma})(x)$ and
$\wh{c} = \apply{Horner}{p_{\pi} \oplus p_{\sigma}, x}$
where $p_{\pi}$ and $p_{\sigma}$ are two polynomials of degree
$n-1$, \Lemma{lemma:HornerSum} yields
\begin{equation*}
|c-\wh{c}|
\leq \gamma_{2n-1} (\wt{p_{\pi}}+\wt{p_{\sigma}})(x)
\leq (1+\uu)^{2n-1} \gamma_{2n-1} \wh{b}.
\end{equation*}
From~(\ref{rel:GammaK}) and~(\ref{rel:StdModel}) it follows that
\begin{equation*}
|c-\wh{c}|
\leq (1+\uu)^{2n} \wh{\gamma}_{2n-1}  \wh{b}
\leq (1+\uu)^{2n+1} \fl(\wh{\gamma}_{2n-1}\wh{b}).
\end{equation*}
Finally we use relation~(\ref{rel:CompBound}) to obtain the error bound.
\end{proof}

\begin{remark}
  \Lemma{lemma:DynBoundErrorTerm} allows us to compute a validated
  error bound for the computed correcting term $\wh{c}$. We apply this
  result twice to derive next \Theorem{prop:DynErrorBounds}. 
  First with \Lemma{lemma:FaithfulRoundingCompHorner} it
  yields the expected dynamic condition for faithful rounding.  
  Then from the EFT for the Horner algorithm (\Theorem{theorem:EFTHorner}) we know
  that $p(x) = \wh{r}+c$. Since $\wb{r} = \wh{r}\oplus\wh{c}$, we deduce
  $|\wb{r} - p(x)| = |(\wh{r}\oplus\wh{c}) - (\wh{r}+\wh{c})+(\wh{c}-c)|$.
  Hence we have
  \begin{equation}
    |\wb{r} - p(x)| \leq |(\wh{r}\oplus\wh{c}) - (\wh{r}+\wh{c})| + |(\wh{c}-c)|.
  \end{equation}
  The first term $|(\wh{r}\oplus\wh{c}) - (\wh{r}+\wh{c})|$ in the previous inequality
  is basically the absolute rounding error that occurs when computing $\wb{r}=\wh{r}\oplus\wh{c}$.
  Using only the bound (\ref{rel:StdModel}) of the standard model of
  floating point arithmetic, it could be bounded by $\uu |\wb{r}|$. But here
  we benefit again from error free transformations using algorithm
  $\algo{TwoSum}$ to compute the actual rounding error exactly, which 
  leads to a sharper error bound. Next \Rel{rel:DynErrorBound}
  improves the dynamic bound presented in \cite{grll05}.
\end{remark}

\begin{theorem} \label{prop:DynErrorBounds}
Consider a polynomial $p$ of degree $n$ with floating point coefficients,
and $x$ a floating point value. Let $\wb{r}$ be the computed value,
$\wb{r}=\apply{CompHorner}{p,x}$ (\Algo{algo:CompHorner}) and let $\wh{\alpha}$ be
the error bound defined by \Rel{rel:DynBoundErrorTerm}.   
\begin{enumerate}[i)]
\item If $\wh{\alpha} < \frac{\uu}{2}|\wb{r}|$, then $\wb{r}$ is a faithful
  rounding of $p(x)$ \label{rel:DynFaithfulBound}.
\item Let $e$ be the floating point value such that $\wb{r}+e = \wh{r}+\wh{c}$, \ie
$[\wb{r}, e] = \apply{TwoSum}{\wh{r},\wh{c}}$, where $\wh{r}$ and
  $\wh{c}$ are defined by \Algo{algo:CompHorner}. The absolute error of
  the computed result $\wb{r} = \apply{CompHorner}{p,x}$ is bounded as follows,
  \begin{equation} \label{rel:DynErrorBound}
    |\wb{r}-p(x)| \leq \fl \left(\frac{\wh{\alpha}+|e|}{1-2u}\right) := \wh{\beta}.
  \end{equation}
\end{enumerate}
\end{theorem}

\begin{proof}
The first proposition follows directly from
\Lemma{lemma:FaithfulRoundingCompHorner}.

By hypothesis $\wb{r}= \wh{r}+\wh{c}-e$, and from
\Theorem{theorem:EFTHorner} we have $p(x) = \wh{r}+c$, thus
\begin{equation*}
|\wb{r}-p(x)| = |\wh{c}-c-e| \leq |\wh{c}-c| + |e| \leq \wh{\alpha} + |e|.
\end{equation*}
From~(\ref{rel:StdModel}) and~(\ref{rel:CompBound}) it follows that
\begin{equation*}
|\wb{r}-p(x)| \leq (1+\uu) \fl(\wh{\alpha} + |e|)
              \leq \fl \left(\frac{\wh{\alpha}+|e|}{1-2u}\right);
\end{equation*}
which proves the second proposition.
\end{proof}

From \Theorem{prop:DynErrorBounds} we deduce the following algorithm.
It computes the compensated result $\wb{r}$ together with the
validated error bound $\wh{\beta}$. Moreover, the boolean value
$\textsf{isfaithful}$ is set to true if and only if the result is proved
to be faithfully rounded.

\begin{algorithm} \label{algo:FaithfulRoundingCompHorner}
Compensated Horner algorithm with check of the faithful rounding
\begin{tabbing}
\quad \=\quad  \kill
function $[\wb{r}, \wh{\beta}, \mbox{\textsf{isfaithful}}] = \apply{CompHornerIsFaithul}{p, x}$ \\[0.1cm]
\> $\left[\wh{r}, p_{\pi}, p_{\sigma}\right] = \apply{EFTHorner}{p, x}$ \\[0.1cm]
\> $\wh{c} = \apply{Horner}{p_{\pi} \oplus p_{\sigma}, x}$\\[0.1cm]
\> $\wh{b} = \apply{Horner}{|p_{\pi}| \oplus |p_{\sigma}|, |x|}$\\[0.1cm]
\> $[\wb{r}, e] = \apply{TwoSum}{\wh{r}, \wh{c}}$\\[0.1cm]
\> $\wh{\alpha} = (\wh{\gamma}_{2n-1} \otimes \wh{b}) \oslash (1 \ominus 2(n+1) \otimes \uu)$\\[0.1cm]
\> $\wh{\beta} = (\wh{\alpha} \oplus |e|) \oslash (1-2 \otimes \uu)$\\[0.1cm]
\> $\mbox{\textsf{isfaithful}} = (\wh{\alpha} < \frac{\uu}{2}|\wb{r}|)$
\end{tabbing}
\end{algorithm}

%%%%%%%%%%%%%%%%%%%%%%%%%%%%%%%%%%%%%%%%%%%%%%%%%%%%%%%%%%%%%%%%%%%%%%%%%%%%%%%%%%
\section{Experimental results}
\label{sec:ExperimentalResults}
%%%%%%%%%%%%%%%%%%%%%%%%%%%%%%%%%%%%%%%%%%%%%%%%%%%%%%%%%%%%%%%%%%%%%%%%%%%%%%%%%%
%%%%%%%%%%%%%%%%%%%%%%%%%%%%%%%%%%%%%%%%%%%%%%%%
%\subsection{Implementation issues}
%%%%%%%%%%%%%%%%%%%%%%%%%%%%%%%%%%%%%%%%%%%%%%%%

We consider polynomials $p$ with floating point coefficients and
floating point entries $x$.  
For presented accuracy tests we use Matlab codes for \algo{CompHorner}
(\Algo{algo:CompHorner}) and \algo{CompHornerIsFaithul}
(\Algo{algo:FaithfulRoundingCompHorner}). These Matlab programs are
presented in Appendix~\ref{sec:MatlabCodes}. 
From these Matlab codes, we see that \algo{CompHorner} requires
$O(21n)$ flop and that \algo{CompHornerIsFaithul} requires $O(26n)$
flop.  

For time performance tests previous algorithms are coded in
C language and several test platforms are described in next
\Tab{tab:MeasTimePerf}.

%%%%%%%%%%%%%%%%%%%%%%%%%%%%%%%%%%%%%%%%%%%%%%%%
\subsection{Accuracy tests}
%%%%%%%%%%%%%%%%%%%%%%%%%%%%%%%%%%%%%%%%%%%%%%%%
We start testing the efficiency of faithful rounding with compensated Horner
algorithm and the dynamic control of faithfulness. 
Then we focus more on both the \apriori and dynamic bounds with two
other test sets. Three cases may occur when the dynamic test
for faithful rounding in \Algo{algo:FaithfulRoundingCompHorner} is
performed.
\begin{enumerate}
\item The computed result is faithfully rounded and this is ensured by
  the dynamic test. Corresponding plots are green in next figures. 
\item The computed result is actually faithfully rounded but the dynamic
  test fails to ensure this property. Corresponding plots are blue.
\item The computed result is not faithfully rounded and plotted in red
  in this case.
\end{enumerate}
Next figures should be observed in color.    

%%%%%%%%%%%%%%%%%%%%%%%%%%%%%%%%
\subsubsection{Faithful rounding with compensated Horner}
%%%%%%%%%%%%%%%%%%%%%%%%%%%%%%%%
\begin{figure}
\begin{center}
\scalebox{0.75}{\includegraphics{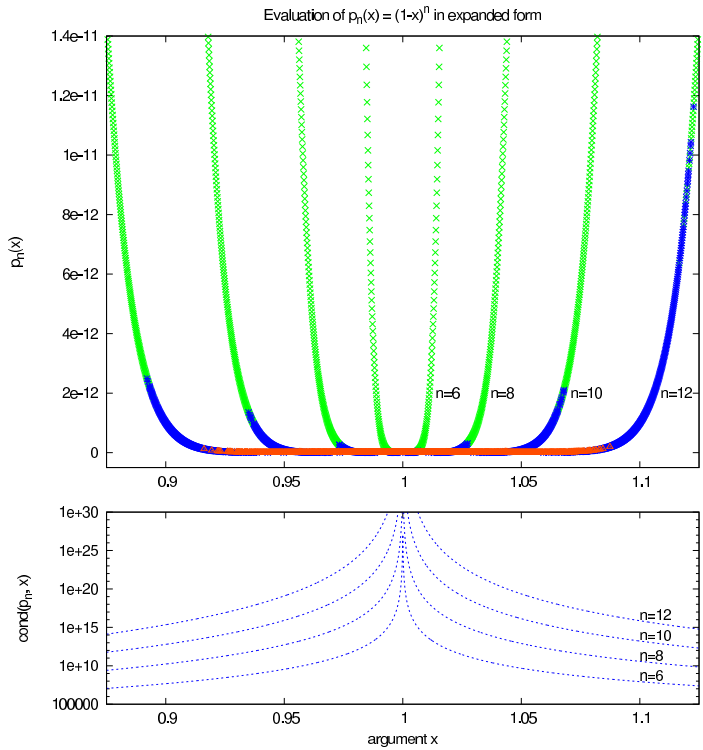}}
\vspace{-0.5cm}
\end{center}
\caption{We report the evaluation of polynomials $p_n$ near
the multiple root $x=1$ with the compensated Horner algorithm
(\algo{CompHornerIsFaithul}) and for
multiplicity $n=6, 8, 10, 12$. Each
evaluation proved to be faithfully rounded thanks to the dynamic test
is reported with a green cross. The faithful evaluations that are not
detected to be so with the dynamic test are represented in blue.
Finally, the evaluations that are not faithfully rounded are reported
in red. The lower frame represents the condition number with respect
to the argument $x$.}
\label{fig:CHFaithBinom}
\end{figure}

In the first experiment set, we evaluate the expanded form of
polynomials $p_n(x) = (1-x)^n$, for degree $n=6,8,10,12$, at $2048$
equally spaced floating point entries being near the multiple root
$x=1$. 
These evaluations are extremely ill-conditioned since 
\begin{equation*}
\cond(p_n, x) = \left| \frac{1+|x|}{1-x} \right|^n.
\end{equation*}
These condition numbers are plotted in the lower frame of
\Fig{fig:CHFaithBinom} while $x$ varies around the root. 
These huge values have a sense since
polynomials $p$ are exact in IEEE-754 double precision. 
Results are reported on \Fig{fig:CHFaithBinom}. 
The well known relation between the lost of accuracy and
the nearness and the multiplicity of the root, \ie the increasing of the
condition number, is clearly illustrated. These results also
illustrate that the dynamic bound becomes more pessimistic as the
condition number increases. In next figures the horizontal axis does
not represent the $x$ entry range anymore but the condition number
which governs the whole behavior.\\  

\begin{figure}
\begin{center}
\scalebox{0.75}{\includegraphics{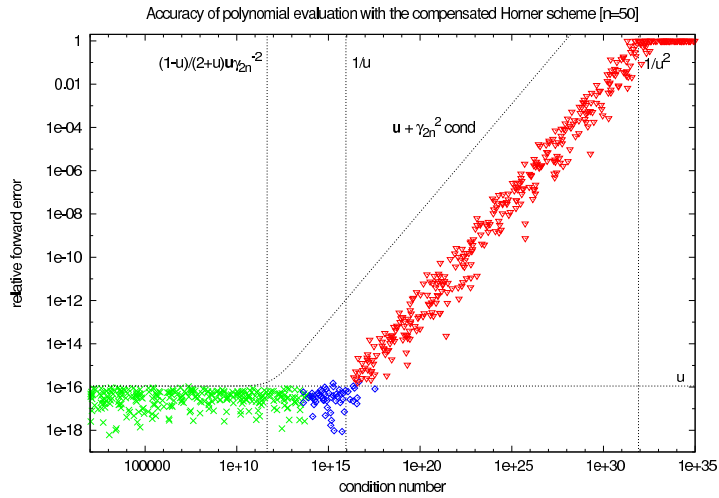}}
\vspace{-0.5cm}
\end{center}
\caption{We report the relative accuracy of every polynomial
  evaluation ($y$ axis) with respect to the condition number ($x$
  axis). Evaluation is performed with
  \algo{CompHornerIsFaithul} (\Algo{algo:FaithfulRoundingCompHorner}). 
  The color code is the same as for \Fig{fig:CHFaithBinom}. 
  Leftmost vertical line is the \apriori sufficient condition
  (\ref{rel:APrioriCriterion}) while the right one marks the inverse of
  the working precision $\uu$. Broken line is the \apriori accuracy bound
  (\ref{rel:REBCompHorner}).}
\label{fig:CHFaithGen}
\end{figure}

For the next experiment set, we first designed a generator of
arbitrary ill-conditioned polynomial evaluations. 
It relies on the condition number definition~(\ref{rel:CondPoly}).
Given a degree $n$, a floating point argument $x$ and a targeted
condition number $C$, it generates a polynomial $p$ with floating
point coefficients such that $\cond(p, x)$ has the same order of
magnitude as $C$. The principle of the generator is the following.
\begin{enumerate}
\item $\lfloor n/2 \rfloor$ coefficients are randomly selected and
  generated such that $\wt{p}(x) = \sum |a_i||x|^i \approx C$,
\item the remaining coefficients are generated ensuring $|p(x)|
  \approx 1$ thanks to high accuracy computation.
\end{enumerate}
Therefore we obtain polynomials $p$ such that
$\cond(p, x) = \wt{p}(x) / |p(x)| \approx C$, for arbitrary values of $C$.

In this test set we consider generated polynomials of degree $50$ whose
condition numbers vary from about $10^2$ to $10^{35}$. These huge
condition numbers again have a sense here since the coefficients and
the argument of every polynomial are floating point numbers.
The results of the tests performed with \algo{CompHornerIsFaithul}
 (\Algo{algo:FaithfulRoundingCompHorner}) are 
reported on \Fig{fig:CHFaithGen}. As expected  every
polynomial with a condition number smaller than the \apriori
bound~(\ref{rel:APrioriCriterion}) is faithfully evaluated with
\Algo{algo:FaithfulRoundingCompHorner} ---green plots at the left of
the leftmost vertical line. 

On \Fig{fig:CHFaithGen} we also see that evaluations with faithful
rounding appear for condition 
numbers larger than the \apriori bound~(\ref{rel:APrioriCriterion}) ---
green and blue plots at the right of the leftmost vertical line. 
As expected a large part of these cases are detected by the dynamic
test introduced in  \Theorem{prop:DynErrorBounds} ---the green ones. 
Next experiment set comes back to this point.
We also notice that the compensated Horner algorithm produces accurate
evaluations for condition numbers up to about $1/\uu$ ---green and
blue plots.

%%%%%%%%%%%%%%%%%%%%%%%%%%%%%%%%
\subsubsection{Significance of the dynamic error bound}
%%%%%%%%%%%%%%%%%%%%%%%%%%%%%%%%

\begin{figure}
\begin{center}
\scalebox{0.8}{\includegraphics{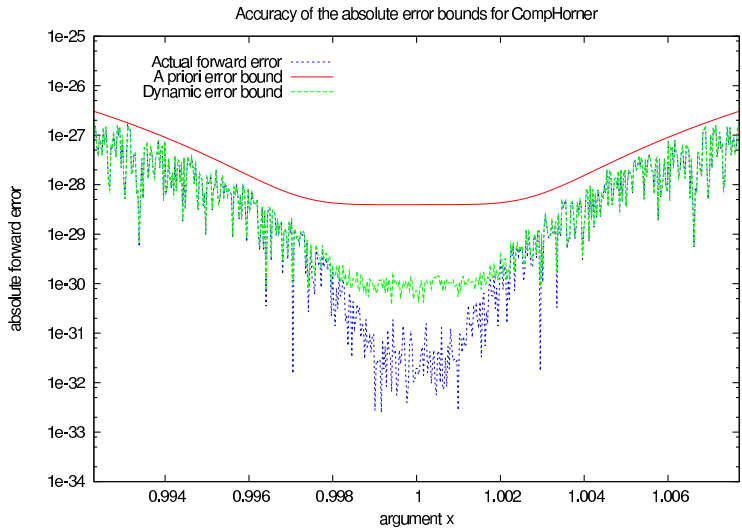}}
\vspace{-0.5cm}
\end{center}
\caption{The dynamic error bound~(\ref{rel:DynErrorBound}) compared
to the \apriori bound~(\ref{rel:EBCompHorner}) and to the actual forward
error ($p(x)=(1-x)^5$ for $400$ entries on the $x$ axis).}
\label{fig:CHBounds}
\end{figure}

We illustrate the significance of the dynamic error bound~(\ref{rel:DynErrorBound}),
compared to the \apriori error bound~(\ref{rel:EBCompHorner}) and to the actual
forward error. We evaluate the expanded form of $p(x) = (1-x)^5$ for $400$
points near $x=1$. For each value of the argument $x$, we compute 
$\apply{CompHorner}{p, x}$ (\Algo{algo:CompHorner}), the associated
dynamic error bound (\ref{rel:DynErrorBound}) and the actual forward
error. The results are reported on \Fig{fig:CHBounds}.

As already noticed, the closer the argument is
to the root $1$ (\ie, the more the condition number increases), the more
pessimistic becomes the \apriori error bound. Nevertheless our dynamic error
bound is more significant than the \apriori error bound  as it takes
into account the rounding errors that occur during the computation.

%%%%%%%%%%%%%%%%%%%%%%%%%%%%%%%%%%%%%%%%%%%%%%%%
\subsection{Time performances}
%%%%%%%%%%%%%%%%%%%%%%%%%%%%%%%%%%%%%%%%%%%%%%%%

\begin{table}
\caption{Measured time performances for \algo{CompHorner}, \algo{CompHornerIsFaithul}
and \algo{DDHorner}. GCC denotes the GNU Compiler Collection and ICC denotes the Intel
C/C++ Compiler.}
\label{tab:MeasTimePerf}
\begin{center}
\renewcommand{\arraystretch}{1.2}
\begin{tabular}{|ll|c|c|c|}
\cline{3-5}
\multicolumn{2}{c|}{}
  & $\frac{\algo{CompHorner}}{\algo{Horner}}$
  & $\frac{\algo{CompHornerIsFaith}}{\algo{Horner}}$
  & $\frac{\algo{DDHorner}}{\algo{Horner}}$ \\ \hline
\small{Pentium 4, 3.00 GHz} & \small{GCC 3.3.5} & 3.77 & 5.52 & 10.00\\
                            & \small{ICC 9.1}   & 3.06 & 5.31 &  8.88\\ \hline
\small{Athlon 64, 2.00 GHz} & \small{GCC 4.0.1} & 3.89 & 4.43 & 10.48\\ \hline
\small{Itanium 2, 1.4 GHz}  & \small{GCC 3.4.6} & 3.64 & 4.59 &  5.50\\
                            & \small{ICC 9.1}   & 1.87 & 2.30 &  8.78\\ \hline
\multicolumn{2}{c|}{}
  & $\sim 2 - 4$
  & $\sim 4 - 6$
  & $\sim 5 - 10$\\ \cline{3-5}
\end{tabular}
\end{center}
\end{table}

All experiments are performed using IEEE-754 double precision. Since the
double-doubles~\cite{LiHB:01,XBLAS:02} are usually considered as the most
efficient portable library to double the IEEE-754 double precision, we consider it
as a reference in the following comparisons. For our purpose, it suffices to know
that a double-double number $a$ is the pair $(a_h, a_l)$ of IEEE-754 floating
point numbers with $a = a_h + a_l$ and $|a_l| \leq \uu |a_h|$. This property
implies a renormalisation step after every arithmetic operation with double-double
values. We denote by \algo{DDHorner} our implementation of the Horner algorithm
with the double-double format, derived from the implementation
proposed in \cite{XBLAS:02}.

We implement the three algorithms \algo{CompHorner}, \algo{CompHornerIsFaith}
and \algo{DDHorner} in a C code to measure their overhead compared to
the \algo{Horner} algorithm. We program these tests straightforwardly with no
other optimization than the ones performed by the compiler. All timings are
done with the cache warmed to minimize the memory traffic over-cost.

We test the running times of these algorithms for different
architectures with different compilers as described in
\Tab{tab:MeasTimePerf}. Our measures are performed with polynomials
whose degree vary from 5 to 200 by step of 5. For each algorithm, 
we measure the ratio of its computing time over the computing time of the classic
Horner algorithm; we display the average time ratio over all test cases in
\Tab{tab:MeasTimePerf}.\\

The results presented in \Tab{tab:MeasTimePerf} show that the slowdown factor
introduced by \algo{CompHorner} compared to the classic \algo{Horner}
roughly varies between 2 and 4. The same slowdown factor varies between 4 and 6 for
\algo{CompHornerIsFaithul} and between 5 and 10 for \algo{DDHorner}. We can see
that \algo{CompHornerIsFaithul} runs a most 2 times slower than \algo{CompHorner}:
the over-cost due to the dynamic test for faithful rounding is therefore quite
reasonable. Anyway \algo{CompHorner} and \algo{CompHornerIsFaithul} run both
significantly faster than \algo{DDHorner}.

\begin{remark}
We provide time ratios for IA'64 architecture (Itanium 2). Tested
algorithms take benefit from IA'64 instructions, \eg \algo{fma}, but
are not described in this paper.   
\end{remark}

%%%%%%%%%%%%%%%%%%%%%%%%%%%%%%%%%%%%%%%%%%%%%%%%%%%%%%%%%%%%%%%%%%%%%%%%%%%%%%%%%%
\section{Conclusion}
\label{sec:Conclusion}
%%%%%%%%%%%%%%%%%%%%%%%%%%%%%%%%%%%%%%%%%%%%%%%%%%%%%%%%%%%%%%%%%%%%%%%%%%%%%%%%%%
Compensated Horner algorithm yields more accurate polynomial
evaluation than the classic Horner iteration. Its accuracy behavior is
similar to an Horner iteration performed in a doubled working
precision. Hence compensated Horner may perform a faithful polynomial
evaluation with IEEE-754 floating point arithmetic in the ``round to
the nearest'' rounding mode. An \apriori sufficient condition with
respect on the condition number that ensures such faithfulness  
 has been defined thanks to the error free transformations. 

These error free transformations also allow us to derive a dynamic
sufficient condition that is more significant to check for faithful
rounding with compensated Horner algorithm. 

It is interesting to remark here that the significance of
this dynamic bound can be improved easily ---how to transform blue
plots in green ones? Whereas bounding the error
in the computation of the (polynomial) correcting term in
\Rel{rel:DynBoundErrorTerm}, a good approximate of the actual error could be
computed (applying again \algo{CompHorner} to the correcting term). Of
course such extra computation will introduce more running time
overhead not necessary useful ---green plots are here! So it suffices
to run such extra (but costly) checking only if the previous
dynamic one fails (a similar strategy as in dynamic filters for
geometric algorithms). 

Compared to the classic Horner algorithm, experimental results exhibit
reasonable over-costs for accurate polynomial evaluation (between 2 and
4) and even for this computation with a dynamic checking for
faithfulness (between 4 and 6). Let us finally remark than such
computation that provides as accuracy as if the working precision is
doubled and a faithfulness checking is no more costly in term 
of running time than the ``double-double'' counterpart without any
check. 

Future work will be to consider subnormals results and
also an adaptative algorithm that ensure faithful rounding
for polynomials with an arbitrary condition number. 

\bibliographystyle{abbrv}
\bibliography{LaLo06}

\begin{thebibliography}{10}

\bibitem{Dekk:71}
T.~J. Dekker.
\newblock A floating-point technique for extending the available precision.
\newblock {\em Numer. Math.}, 18:224--242, 1971.

\bibitem{grll05}
S.~Graillat, P.~Langlois, and N.~Louvet.
\newblock Compensated {H}orner scheme.
\newblock Technical report, University of Perpignan, France, July 2005.

\bibitem{LiHB:01}
Y.~Hida, X.~S. Li, and D.~H. Bailey.
\newblock Algorithms for quad-double precision floating point arithmetic.
\newblock In N.~Burgess and L.~Ciminiera, editors, {\em Proceedings of the 15th
  Symposium on Computer Arithmetic, Vail, Colorado}, pages 155--162, Los
  Alamitos, CA, USA, 2001. Institute of Electrical and Electronics Engineers.

\bibitem{ASNA:02}
N.~J. Higham.
\newblock {\em Accuracy and Stability of Numerical Algorithms}.
\newblock Society for Industrial and Applied Mathematics, Philadelphia, PA,
  USA, second edition, 2002.

\bibitem{IEEE:85}
{IEEE Standards Committee 754}.
\newblock {\em {IEEE} Standard for binary floating-point arithmetic,
  {ANSI}/{IEEE} {S}tandard 754-1985}.
\newblock Institute of Electrical and Electronics Engineers, Los Alamitos, CA,
  USA, 1985.
\newblock Reprinted in {SIGPLAN} {N}otices, 22(2):9-25, 1987.

\bibitem{Knut:98}
D.~E. Knuth.
\newblock {\em The Art of Computer Programming: Seminumerical Algorithms},
  volume~2.
\newblock Ad{\-d}i{\-s}on-Wes{\-l}ey, Reading, MA, USA, third edition, 1998.

\bibitem{XBLAS:02}
X.~S. Li, J.~W. Demmel, D.~H. Bailey, G.~Henry, Y.~Hida, J.~Iskandar, W.~Kahan,
  S.~Y. Kang, A.~Kapur, M.~C. Martin, B.~J. Thompson, T.~Tung, and D.~J. Yoo.
\newblock Design, implementation and testing of extended and mixed precision
  {BLAS}.
\newblock {\em ACM Trans. Math. Software}, 28(2):152--205, 2002.

\bibitem{OgRO:05}
T.~Ogita, S.~M. Rump, and S.~Oishi.
\newblock Accurate sum and dot product.
\newblock {\em SIAM J. Sci. Comput.}, 26(6):1955--1988, 2005.

\bibitem{Prie:91}
D.~M. Priest.
\newblock Algorithms for arbitrary precision floating point arithmetic.
\newblock In P.~Kornerup and D.~W. Matula, editors, {\em Proceedings of the
  10th {IEEE} Symposium on Computer Arithmetic (Arith-10),Grenoble, France},
  pages 132--144, Los Alamitos, CA, USA, 1991. Institute of Electrical and
  Electronics Engineers.

\bibitem{RuOO06}
S.~M. Rump, T.~Ogita, and S.~Oishi.
\newblock Accurate summation.
\newblock Technical report, Hamburg University of Technology, Germany, Nov.
  2005.

\end{thebibliography}

\section{Appendix}
\label{sec:MatlabCodes}
%%%%%%%%%%%%%%%%%%%%%%%%%%
Accuracy tests use next Matlab codes for algorithms \Algo{algo:CompHorner}
(\algo{CompHorner}) and \Algo{algo:FaithfulRoundingCompHorner}
(\algo{CompHornerIsFaithul}). Following Matlab convention,
$p$ is represented as a vector \texttt{p} such that
$p(x) = \sum_{i=0}^{n} \mathtt{p}(n-i+1) x^i$.
We also recall that Matlab \texttt{eps} denotes the machine epsilon, which is
the spacing between $1$ and the next larger floating point number, hence
$\uu = \mathtt{eps}/2$.

\begin{minipage}{\linewidth}
\vspace{0.3cm}
%%%%%%%%%%%%%%%%%%%%%%%%%%%%%%
\begin{minipage}[t]{0.47\linewidth}
\begin{algorithm} \label{algo:CodeCompHorner}
Code for \Algo{algo:CompHorner}.
\begin{small}
\let\ttdefault\sfdefault
\begin{verbatim}
function r = CompHorner(p, x)
  n = length(p)-1; % degree of p
  [xh, xl] = Split(x);
  r = p(1); c = 0.0;    
  for i=2:n+1
    %[r, pi] = TwoProd(r, x)
    p = r*x;
    [rh, rl] = Split(r);
    pi = rl*xl-(((p-rh*xl)-rl*xh)-rh*xl);
    %[r, sigma] = TwoSum(r, p(i))
    r = p+p(i);
    t = r-p;
    sigma = (p-(r-t))+(p(i)-t);
    % Computation of the correcting term
    c = c*x+(pi+sig);
  end
  % Final correction of the result
  r = r+c;
\end{verbatim}

\end{small}
\end{algorithm}
\end{minipage}
%%%%%%%%%%%%%%%%%%%%%%%%%%%%%%
\begin{minipage}[t]{0.47\linewidth}
\begin{algorithm} \label{algo:CodeFaithfulRoundingCompHorner}
Code for \Algo{algo:FaithfulRoundingCompHorner}.
\begin{small}
\let\ttdefault\sfdefault
\begin{verbatim}
function [r, beta, isfaith] = CompHornerIsFaithul(p, x)
  n = length(p)-1; % degree of p
  [xh, xl] = Split(x);
  absx = abs(x);
  r = p(1); c = 0.0; beta = 0.0;
  for i=2:n+1
    % [r, pi] = TwoProd(r, x)
    p = r*x;
    % [rh, rl] = Split(r);
    pi = rl*xl-(((p-rh*xl)-rl*xh)-rh*xl);
    % [r, sigma] = TwoSum(r, p(i))
    r = p+p(i);
    t = r-p;
    sigma = (p-(r-t))+(p(i)-t);
    % Computation of the correcting term
    c = c*x+(pi+sig);
    b = b*absx+(abs(pi)+ abs(sig));
  end
  % Final correction of the result
  [r, e] = TwoSum(r,c);
  % Check for faithful rounding
  alpha = gam(2*n-1)*b / (1-(n+1)*eps);
  isfaith = alpha < 0.25*eps*abs(r);
  % Absolute error bound
  beta = (alpha + abs(e))/(1-2*u);
\end{verbatim}
\end{small}
\end{algorithm}
\end{minipage}
%%%%%%%%%%%%%%%%%%%%%%%%%%%%%%
\vspace{0.3cm}
\end{minipage}

\end{document}